\documentclass[10pt,reqno]{amsart}
\usepackage{amssymb,accents,calrsfs}
\def\noi{\noindent}

\newtheorem{Thm}{Theorem}[section]
\newtheorem{Def}[Thm]{Definition}
\newtheorem{Lm}[Thm]{Lemma}
\newtheorem{Prop}[Thm]{Proposition}
\newtheorem{Cor}[Thm]{Corollary}

\def\cal{\mathcal}
\def\Bbb{\mathbb}
\def\mf{\mathfrak}

\def\<{\langle}
\def\>{\rangle}

\def\a{\alpha}
\def\b{\beta}
\def\d{\delta}

\def\th{\theta}

\def\l{\lambda}
\def\L{\Lambda}

\def\Re{\Bbb R}
\def\F{\Bbb F}
\def\C{\Bbb C}
\def\Z{\Bbb Z}

\def\X{\mf X}
\def\T{\mf T}
\def\R{\mf R}

\def\O{\mf O}
\def\K{\mf K}
\def\A{\cal A}

\def\G{\mf h}

\def\W{\ring W}
\def\w{\ring w}
\def\RR{\ring R}
\def\Q{\ring Q}

\def\Gc{\ring {\mf h}}

\def\p{\mf p}
\def\t{\mf t}

\begin{document}
\title[]
{Nonsymmetric Macdonald polynomials and matrix coefficients for
unramified principal series}
\author{Bogdan Ion}

\thanks{Department of Mathematics, University of Michigan, Ann
Arbor, MI 48109.}
\thanks{E-mail address: {\tt bogdion@umich.edu}}\thanks{Supported in part by a Rackham Faculty
Research Fellowship}
\begin{abstract}
We show how a certain limit of the nonsymmetric Macdonald
polynomials appears in the representation theory of semisimple
groups over $\mf p$--adic fields  as matrix coefficients for the
unramified principal series representations. The result is the
nonsymmetric counterpart of a classical result relating the same
limit of the symmetric Macdonald polynomials to zonal spherical
functions on groups of $\mf p$--adic type.
\end{abstract}
\maketitle

\thispagestyle{empty}
\section*{Introduction}

The theory of zonal spherical functions for semisimple groups of
$\mf p$--adic type was completed by the early 1970's through the
work of Satake \cite{satake} and Macdonald \cite{macdonald} (see
also \cite{casselman}). Unlike the case of real semisimple groups,
an explicit formula for the values of the zonal spherical
functions for the groups of $\mf p$--adic type was obtained and
the Plancherel measure on the set of positive definite spherical
functions (or, equivalently, on the set of irreducible unitary
spherical representations) was computed. Later on, partially
inspired by this formula, Macdonald constructed a remarkable
family of orthogonal polynomials $P_\l(q,t)$ associated to any
finite, irreducible root system $\RR$. These polynomials are Weyl
groups invariant functions which depend rationally on two sets of
parameters $q$ and $t$ and are indexed by the anti--dominant
elements of the weight lattice of $\RR$. Their limit
$P_\l(\infty,t)$, as $q$ approaches infinity, coincides up to a
scalar factor with the previously obtained formula for the
spherical functions associated to a group of $\mf p$--adic type
carefully chosen to match (in a very precise sense) the root
system $\RR$.

In a further development, Opdam (in the so--called differential
setting), Macdonald (in the difference $q,t$--setting), Cherednik
(general case) and Sahi (for nonreduced root systems) introduced
another remarkable family of orthogonal polynomials $E_\l(q,t)$
also associated to finite root systems but indexed now by the full
weight lattice. In contrast with $P_\l(q,t)$, the polynomials
$E_\l(q,t)$ are nonsymmetric (i.e. not invariant under the action
of the Weyl group). The nonsymmetric Macdonald polynomials became
increasingly important from various points of view. On one hand,
they enjoy many of the properties of the symmetric Macdonald
polynomials. For example they satisfy versions of Macdonald's norm
and evaluation--duality conjectures and they provide
simplifications in the proofs of the original conjectures (see,
for example, \cite{macbook}). On the other hand, they were the key
ingredient in many recent developments in the theory of orthogonal
polynomials and related combinatorics (for example, in the proof
of Macdonald's positivity conjecture for Jack polynomials by Knop
and Sahi), and for the harmonic analysis and representation theory
of double affine Hecke algebras (through the work of Cherednik).

Despite the growing importance of nonsymmetric Macdonald
polynomials not much it is known with respect to their
interpretation in classical representation--theoretical terms. The
main difficulty in integrating them within the classical framework
lies, on one hand, on our lack of explicit formulas for
nonsymmetric Macdonald polynomials and, on the other hand, on the
fact that most special functions (Weyl characters, spherical
functions, conformal blocks) appearing in representation--theory
are Weyl group invariant. However, a first result in this
direction was obtained in \cite{ion2} were it was shown that in
the limit $E_\l(q,\infty)$, as $t$ approaches infinity, the
nonsymmetric Macdonald polynomials coincide with Demazure
characters for basic representations of affine Kac--Moody groups,
and in the limit $E_\l(\infty,\infty)$, as both $q$ and $t$
approach infinity, the nonsymmetric Macdonald polynomials coincide
with Demazure characters for irreducible representations of
compact Lie groups.

The goal of this paper is to give a representation--theoretical
interpretation for the limit $E_\l(\infty,t)$, as $q$ tends to
infinity, of the nonsymmetric Macdonald polynomials. As mentioned
above, in this limit, their symmetric counterparts
$P_\l(\infty,t)$ describe the values of spherical functions on
groups $G$ over $\mf p$--adic fields, with respect to $K$, a
special maximal compact subgroup of $G$. Keeping in mind that the
zonal spherical functions are in fact matrix coefficients
(corresponding to pairs of $K$--fixed vectors) for irreducible
spherical representations, our main results, Theorem \ref{thm2}
and Theorem \ref{thm1}, attach a similar meaning to the
polynomials $E_\l(\infty,t)$: they give the values of matrix
coefficients for unramified principal series representations of
$G$. The two vectors involved in the computation of each matrix
coefficient are: one, a (essentially unique) $K$--fixed vector,
and the other, a carefully chosen Iwahori fixed vector (see
Section \ref{coeff} for details).

In contrast with the corresponding result for symmetric Macdonald
polynomials where the identification was established by comparing
the explicit formulas on both sides, in our case, the result is
made possible by  the recursion formula for nonsymmetric Macdonald
polynomials obtained in \cite{ion}. With the exception of Sections
1 and 3 which form a brief exposition of the results from the
theory of Macdonald polynomials which will be used later on, the
rest of the paper is concerned with establishing the framework
within the theory of groups over $\mf p$--adic fields which will
make possible the identification between nonsymmetric Macdonald
polynomials and certain matrix coefficients. Section 2 is again
expository, presenting the main structural results on $\mf
p$--adic groups as they follow from the Bruhat--Tits theory.
Sections 4 and 5, which contain the main technical results, are
inspired by the (unpublished) work of Bernstein on projective
generators for blocks (in our case unramified principal series
representations) of the category of smooth representations of $\mf
p$--adic groups. The idea of realizing the Iwahori--Hecke algebras
as endomorphism algebras of such a projective generator appeared
initially (for split groups) in \cite{chriss}. Since we work here
with non necessarily split groups and we need a slightly different
variant of their result we choose to rework here some of the
arguments. Section 6 concludes the proof of our main result.

\section{Affine root systems and Weyl groups}\label{affine}

\subsection{} Let $\RR\subset
\Gc^*$ be a finite, irreducible, not necessarily reduced, root
system of rank $n$, and let $\RR^\vee\subset \Gc$ be the dual root
system. We denote by $\{\a_i\}_{1\leq i\leq n}$ a basis of $\RR$
(whose elements will be called simple roots); the corresponding
elements $\{\a_i^\vee\}_{1\leq i\leq n}$ of $\RR^\vee$ will be
called simple coroots. If the root system is nonreduced, let us
arrange that $\a_n$ is the unique simple root such that $2\a_n$ is
also a root. The root $\th$  is defined as the {\sl highest short
root} in $\RR$ if the root system is reduced, or as the {\sl
highest root} if the root system is nonreduced.

Our choice of basis determines a subset $\RR^+$ of $\RR$ whose
elements are called positive roots; with the notation
$\RR^-:=\RR^+$ we have $\RR=\RR^+\cup \RR^-$. As usual, $\Q
=\oplus_{i=1}^n\Z\a_i$ denotes the root lattice of $\RR$. Let
$\{\l_i\}_{1\leq i\leq n}$ and $\{\l_i^\vee\}_{1\leq i\leq n}$ be
the fundamental weights, respectively the fundamental coweights
associated to $\RR^+$, and denote by $P=\oplus_{i=1}^n\Z\l_i$ the
weight lattice. An element of $P$ will be called dominant if it is
a linear combination of the fundamental weights with non--negative
integer coefficients. Similarly an anti--dominant weight is a
linear combination of the fundamental weights with non--positive
integer coefficients.

The real vector space ${\Gc}^*$ has a canonical scalar product
$(\cdot, \cdot)$ which we normalize such that it gives square
length 2 to the short roots in $\RR$ (if there is only one root
length we consider all roots to be short); if $\RR$ is not reduced
 we normalize the scalar product such that the roots have square
length 1, 2 or 4. We will use $\RR_s$ and $\RR_\ell$ to refer to
the short and respectively long roots in $\RR$; if the root system
is nonreduced we will also use $\RR_m$ to refer to the roots of
length 2. We will identify the vector space $\Gc$ with its dual
using this scalar product. Under this identification
$\a^\vee=2\a/(\a,\a)$ for any root $\a$.

To any finite root system as above we will associate an {\sl
affine root system} $R$. Let ${\rm Aff}(\Gc)$ be the space of
affine linear transformations of $\Gc$. As a vector space, it can
be identified to $\Gc^*\oplus \Re\d$ via
$$
(f+c\d)(x)=f(x)+c, \quad \text{for~} f\in \Gc^*, ~x\in\Gc
\text{~and~} c\in \Re
$$
Assume first that $\RR$ is reduced, and let $r$ denote the maximal
number of laces connecting two vertices in its Dynkin diagram.
Then,
$$R:=(\RR_s+\Z\d)\cup(\RR_\ell +r\Z\d)\subset \Gc^*\oplus \Re\d$$
If the finite root system $\RR$ is nonreduced then
$$R:=(\RR_s+\frac{1}{2}\Z\d)\cup(\RR_m+\Z\d)\cup(\RR_\ell +\Z\d)$$

The set of affine positive roots $R^+$ consists of affine roots of
the form $\a+k\d$ such that $k$ is non--negative if $\a$ is a
positive root, and $k$ is strictly positive if $\a$ is a negative
root. The affine simple roots are $\{\a_i\}_{0\leq i\leq n}$ where
we set $\a_0:=\d-\th$ if $\RR$ is reduced and
$\a_0:=\frac{1}{2}(\d-\th)$ otherwise. In fact, to make our
formulas uniform we set $\a_0:=c_0^{-1}(\d-\th)$, where $c_0$
equals $1$ or 2 depending on whether $\RR$ is reduced or not. The
root lattice of $R$ is defined as $Q=\oplus_{i=0}^n \Z\a_i$.

Abstractly, an affine root system is a subset $\Phi_{\rm
af}\subset {\rm Aff}(V)$ of the space of affine--linear functions
on a real vector vector space $V$, consisting of non--constant
functions which satisfy the usual axioms for root systems. As in
the case of finite root systems, a classification of the
irreducible affine root systems is available (see, for example,
\cite[Section 1.3]{macbook}). The affine root systems $R$ which we
defined above are just a subset of all the irreducible affine root
systems. However, the configuration of vanishing hyperplanes of
elements of an irreducible affine root system $\Phi$ coincides
with the corresponding configuration of hyperplanes  associated to
a {\sl unique} affine root system $R$ as above. Moreover, the
nonreduced affine root systems we consider above contain as
subsystems all the other nonreduced irreducible affine root
systems and also all reduced irreducible affine root systems of
classical type.


\subsection{}

 The scalar product on $\Gc^*$ can be extended to a
non--degenerate bilinear form on the real vector space
$$\G^*:=\Gc^*\oplus \Re\d\oplus \Re\L_0$$ by requiring that $(\d,\Gc^*\oplus \Re\d)=
(\L_0,\Gc^*\oplus \Re\L_0)=0$ and $(\d,\L_0)=1$. Given $\a\in R$
and $x\in \G^*$ let
$$
r_\a(x):=x-\frac{2(x,\a)}{(\a,\a)}\a\
$$
The {\sl affine Weyl group} $W$ is the subgroup of ${\rm
GL}(\G^*)$ generated by all $r_\a$ (the simple reflections
$r_i=r_{\a_i}$ are enough). The {\sl finite Weyl group} $\W$ is
the subgroup generated by $r_1,\dots,r_n$. The bilinear form on
$\G^*$ is equivariant with respect to the affine Weyl group
action.

The affine Weyl group could also be presented as a semidirect
product in the following way: it is the semidirect product of $\W$
and the lattice $\Q$ (regarded as an abelian group with elements
$t_\mu$, where $\mu$ is in $\Q$), the finite Weyl group acting on
the root lattice as follows
$$
\w t_\mu \w^{-1}=t_{\w(\mu)}
$$
Since the finite Weyl group  also acts on the weight lattice, we
can also consider the {\sl extended Weyl group} $W^e$ defined as
the semidirect product between $\W$ and $P$. Unlike the the affine
Weyl group, $W^e$ is not a Coxeter group. However, $W$ is a normal
subgroup of $W^e$ and the quotient is finite. For $\l$ in $P$, the
action of the elements $t_\l$ on $\G^*$ is described below.

For $s$ a real number, $\G^*_s=\{ x\in\G\ ;\ (x,\d)=s\}$ is the
level $s$ of $\G^*$. We have
$$
\G^*_s=\G^*_0+s\L_0=\Gc^*+{\Bbb R}\d+s\L_0\ .
$$
 The action of $W$
preserves each of the $\G^*_s$ and  we can identify each of the
$\G^*_s$ canonically with $\G^*_0$ and obtain an (affine) action
of $W$ on $\G^*_0$, which we can then restrict to $\Gc^*$. By
$s_i\cdot$ and $\tau_\mu\cdot $ we denote the affine action of $W$
on $\Gc^*$
\begin{eqnarray*}
s_0\cdot x   & = & s_\th(x)+c_0^{-1}\th\ ,\\
\tau_\mu\cdot x & = & x  + \mu  \ ,
\end{eqnarray*}

If we examine the orbits of the affine Weyl group $W$ for its
level zero action on the affine root system $R$ we find that if
$\RR$ is reduced there are precisely as many orbits as root
lengths. If $\RR$ is nonreduced and of rank at least two, then the
above action has five orbits
$$
W(2\a_0)=\RR_\ell+2\Z\d+\d, \quad
W(\a_0)=\RR_s+\Z\d+\frac{1}{2}\d,\quad W(\a_1)=\RR_m+\Z\d
$$
$$
 W(2\a_n)=\RR_\ell+2\Z\d\quad\text{and}\quad W(\a_n)=\RR_s+\Z\d
$$
If $\RR$ is nonreduced and of rank one then $\RR_m$ is empty and
therefore the level zero action of the affine Weyl group on $R$
has only four orbits.
\subsection{} Let us introduce a field $\F$ (of parameters) as follows. Let
$t=(t_\a)_{\a\in R}$ be a set of parameters which is indexed by
the set of affine roots and has the property that $t_\a=t_\b$ if
and only if the affine roots $\a$ and $\b$ belong to the same
orbit under the action of the affine Weyl group on the affine root
system. It will be convenient to also have the following
convention: if $\a$ is not an affine root then $t_\a=1$. Let $q$
be another parameter and let $m$ be the lowest common denominator
of the rational numbers $\{(\a_j,\l_k)\ |\ 1\leq j,k\leq n \}$.
The field $\F=\F_{q,t}$ is defined as the field of rational
functions in $q^{\frac{1}{m}}$ and
$t^{\frac{1}{2}}=(t_\a^{\frac{1}{2}})_{\a\in R}$. We will also use
the field of rational functions in
$t^{\frac{1}{2}}=(t_\a^{\frac{1}{2}})_{\a\in R}$ denoted by
$\F_t$. The algebra $\cal R_{q,t}=\F[e^\l;\l\in P]$ is the group
$\F$-algebra of the lattice $P$. Similarly, the algebra $\cal
R_t=\F_t[e^\l;\l\in P]$ is the group $\F_t$-algebra of the lattice
$P$.

If the root system $R$ is reduced then  there are as many distinct
parameters $t_\a$ as root lengths (i.e. at most two). In this
case, for any affine simple root $\a_i$ we will use the notation
$t_i$ to refer to the parameter $t_{\a_i}$. To keep this notation
consistent with the one for nonreduced root systems (detailed
below), we also introduce $t_{01}=t_{02}=t_{03}:=t_0$. If $R$ is
nonreduced then the action of the affine Weyl group on the affine
root system has five orbits $W(2\a_0),~W(\a_0),~W(a_n), ~W(2\a_n)$
and $W(\a_1)$ (note that the last orbit is empty if $R$ has rank
one) and we denote the corresponding parameters by $t_{01},~
t_{02},~ t_{03},~ t_n$ and $t_1=\cdots=t_{n-1}$, respectively.

\section{Groups over $\mf p$--adic fields}\label{groups}
\subsection{} Let $\K$ be a complete, non-archimedean local field with finite
residue field. We denote by $\O$ its ring of integers and by
$\p=(\varpi)\subset \O$ the unique prime ideal of $\O$. The
residue field $\mf f:=\O/\p$ has finite cardinality, denoted by
$\t$. Every element of $\K^\times$ can be uniquely written as
$\varpi^l$ for some integer $l$ and some $u\in \O^\times$, a unit
of  $\O$. The valuation $v:\K\to \Z$ is defined as
$v(\varpi^lu)=l$ (by convention $v(0)=+\infty$). The absolute
value of an element $x$ of $\K$ is defined as $|x|=\t^{-v(x)}$.
The metric induced by the absolute value makes $\K$ into a locally
compact, totally disconnected topological space.


\subsection{} Let ${\bf G}$ be a
connected, $\K$--simple, linear algebraic group defined over $\K$.
We will also assume that $\bf G$ is $\K$--isotropic of $\K$--rank
$n$.   We denote by $G$ the  $\K$--rational points of ${\bf G}$
and the same type of notation will be used for all the linear
algebraic groups defined in this paragraph (they are all defined
over $\K$). Let ${\bf P}$ be a minimal $\K$--parabolic subgroup of
${\bf G}$ and ${\bf S}$ a maximal $\K$--split torus of $\bf G$
contained in ${\bf P}$. The parabolic ${\bf P}$ has the Levi
decomposition ${\bf P=MU}$, where ${\bf M}=Z_{\bf G}({\bf S})$ is
the centralizer of ${\bf S}$ in ${\bf G}$, and ${\bf U}$ is the
unipotent radical of ${\bf P}$. The unique $\K$--parabolic
subgroup of ${\bf G}$ which is opposed to ${\bf P}$ with respect
to ${\bf M}$ is denoted by ${\bf P^-}$, and ${\bf U^-}$ denotes
its unipotent radical. Also, let ${\bf N}=N_{\bf G}({\bf S})$ be
the normalizer of ${\bf S}$.


\subsection{} The finitely generated free $\Z$--modules
$X^*(S):={\rm Hom}_\K(S,\K^\times)$ and $X_*(S):={\rm
Hom}_\K(\K^\times, S)$ are the character group and, respectively,
the cocharacter group of $S$. Since $S$ is split $X_*(S)$ could be
also identified with ${\rm Hom}_\Z(X^*(S), \Z)$ via the perfect
pairing $\<\cdot,\cdot\>: X_*(S)\times X^*(S)\to \Z$ which by
definition satisfies $\psi(\varphi(s))=s^{\<\varphi,\psi\>}$, for
all $s\in \K^\times$, $\varphi\in X_*(S)$ and $\psi\in X^*(S)$.
Consider the real vector space $V:=X_*(S)\otimes_\Z \Re$ and its
dual $V^*:=X^*(S)\otimes_\Z \Re$. With this notation, the
$\Re$--linear extension of the pairing $\<\cdot,\cdot\>$ becomes
the natural pairing between $V$ and $V^*$.

Similarly, consider $X^*(M)={\rm Hom}_\K(M,\K^\times)$ the group
of $\K$--rational characters of $M$, and $X_*(M)={\rm
Hom}_\Z(X^*(M),\Z)$. Since $X^*(M)$ is a full rank sublattice of
$X^*(S)$, the natural pairing between $X_*(M)$ and $X^*(M)$ is
also $\<\cdot,\cdot\>$. We have the following group morphism
$$
v_M: M\to X_*(M)
$$
defined by $\<v_M(m),\psi\>:=v(\psi(m))$, for any $m$ in $M$ and
any $\psi$ in $X^*(M)$. The image, respectively the kernel of
$v_M$ will be denoted by $\L$ and $\ring {}M$, respectively. The
subgroup $\ring {}M$ is a maximal compact and open subgroup of $M$
and, similarly, $\ring {}S:=S\cap \ring {}M$ is a maximal compact
and open subgroup of $S$.

 The inclusion of $S$ into $M$ induces the inclusions
$X_*(S)\subset \L\subset X_*(M)$ which in general are strict. If
${\bf G}$ is $\K$--split, both inclusions are equalities, and if
${\bf G}$ splits over an unramified extension of $\K$ then $\L$
coincides with $X_*(S)$. However, $X_*(M)\otimes_\Z
\Re=\L\otimes_\Z \Re=V$ and the three sets considered above are
lattices of $V$ of rank $n$.

\subsection{}
A character $\chi:M\to \C^\times$ is called unramified if it is
trivial on ${\ring {}M}$. Therefore $X_{\rm nr}(M)$, the set of
unramified characters of $M$, is in bijection with the set of
group morphisms ${\rm Hom}(\L, \C^\times)$ or, equivalently, with
the set of $\C$--algebra morphisms ${\rm Hom_{\C-alg}}(\R,\C)$,
where we denoted by ${\R}$ the group $\C$--algebra of the lattice
$\L$. When we consider an unramified character $\chi$ of $M$, we
will use the same symbol to refer to the corresponding morphisms
in any of the above contexts. We use the notation $\C_\chi$ to
refer to $\C$ with the $\R$--module structure given by the
unramified character $\chi:\R\to \C$.

The lattice $\L$ can be considered both as a subset of $V$ and a
quotient of $M$ and the group structure will be denoted additively
or multiplicatively depending on the situation. For example, if
$\l$ and $\mu$ are elements of $\L\subset V$, then the
corresponding elements of $M$ (determined up to multiplication by
elements of $~\ring {}M$) are denoted by $t_\l$ and $t_\mu$, their
product is $t_\l t_\mu=t_{\l+\mu}$ and $v_M(t_\l)=\l$. For an
element $\l$ of $\L$ we will denote by $\mf e^\l$ the
corresponding element in $\R$. The group structure  of $\R$ will
be also denoted multiplicatively.


\subsection{} Consider ${\cal B}$, the Bruhat--Tits building
associated to $G$, and ${\A}$ the apartment in ${\cal B}$
corresponding to the torus $S$. The apartment $\A$ is an
$n$--dimensional real affine space on which $V$ acts by
translations. We also fix once and for all a special vertex $x_0$
in $\A$. The affine root system $\Phi_{\rm af}=\Phi_{\rm af}({\bf
G},{\bf S}, \K)$ consists of affine linear functions on $\A$ whose
vector parts (which are linear functions on $V$) are elements of
$\Phi=\Phi({\bf G},{\bf S}, \K)\subset X^*(S)$, the $\K$--relative
root system of ${\bf G}$ with respect to the torus $\bf S$. The
vanishing hyperplanes of the affine roots will be referred to as
the {\sl affine hyperplanes in $\A$}. The minimal parabolic $P$
determines a basis $\{a_1,\cdots,a_n\}$ of $\Phi$ and,
consequently, a set of positive roots $\Phi^+$.  We will identify
the linear functions on $V$ (in particular the elements of $\Phi$)
with the affine linear functions on $\A$ which vanish at $x_0$.
Note that $\Phi$ may be a nonreduced root system; in such a case
we arrange that $a_n$ is the unique simple root which belongs to
$\frac{1}{2}\Phi$. Also, we denote by $a_0$ the unique
non--divisible affine root which vanishes on the wall of $\cal C$
not containing $x_0$. The elements of the set $\{a_0,\cdots,a_n\}$
are called the affine simple roots.

Among the connected components of the complement in $\A$ of the
union of affine hyperplanes passing through $x_0$ there is a
unique one, denoted by $\mf C$, with the property that $\mf C\cap
u\cdot \mf C$ contains a translate of $\mf C$ for any $u$ in $U$.
There exists a unique chamber $\cal C$ which is included in $\mf
C$ and has $x_0$ as one of its vertices. The cone $\mf C$ can be
alternatively described as the set of points in $\A$ on which the
positive roots $\Phi^+$ take non--negative values.

The stabilizer of $x_0$ in $G$, denoted by $K$, is called a
special, good, maximal compact subgroup of $G$. The interior
points of $\cal C$ all have the same stabilizer $I$ in $G$, called
the Iwahori subgroup of $G$ attached to $\cal C$.

To any affine root $a\in \Phi_{\rm af}$, Tits associates a
strictly positive integer $d(a)$ which depends on the valuation
data on the root subgroups of $G$. For any of the affine simple
roots $a_i$ define the positive integer $d_i:=d(a_i)+d(2a_i)$ (by
convention $d(2a_i)=0$ if $2a_i$ is not an affine root). For
example, if $\bf G$ is $\K$--split, all the integers $d_i$ are 1.
The integer $d_i$ can be also described in geometric terms as
follows: the wall of $\cal C$ on which $a_i$ vanishes is contained
in the closure of precisely ${\mf t}^{d_i}+1$ chambers of $\cal
B$. To simplify later formulas we will use the notation $\mf t_i$
to refer to ${\mf t}^{d_i}$ and $\mf t_i^*$ to refer to $\mf
t^{d(a_i)-d(2a_i)}$.

\subsection{} Any element $g$ of $N$ stabilizes the apartment $\A$
and it fixes every point of $\A$ if and only if it belongs to
$\ring {}M$. In consequence, we identify the {\sl extended affine
Weyl group} $W^e:=N/{~\ring {}M}$ with a group of affine
transformations of $\A$. The {\sl affine Weyl group} $W$ is the
group generated by the reflections $r_a$ in $\A$, with $a$ an
affine root. The simple reflections $r_i:=r_{a_i}$, $0\leq i\leq
n$, in the walls of $\cal C$ are sufficient to generate $W$. The
affine Weyl group $W$ is a normal subgroup of $W^e$. If the finite
group $\Omega$ is the normalizer in $W^e$ of the chamber $\cal C$,
we have that $W^e=\Omega\ltimes W$. The group $\Omega$ can be
equivalently described as $N_G(I)/I$; if $\bf G$ is
simply--connected, then $\Omega$ is trivial and $W^e=W$.

The subgroup of $W$ generated by the simple reflections
$r_1,\cdots,r_n$ is the {finite Weyl group} $\W=W(\Phi)$. Here we
identified the affine space $\A$ with $V$ by choosing $x_0$ as the
origin. Since $G$ is connected $\W$ is isomorphic to $N/M$. Any
element $\w$ of the finite Weyl group $\W\cong N/M$ has a
representative in $K\cap N$ (which we denote by the same symbol).
Since $~\ring {}M=K\cap M$, the group $\W$ is isomorphic to
$N_K(S)/{~\ring {}M}$ which is a subgroup of $W^e$. We also regard
$\L\cong M/{~\ring {}M}$ as a subgroup of $W^e$. Furthermore, the
extended affine Weyl group $W^e$ can be presented as $\W\ltimes
\L$. From now on we will assume that we have chosen a
representative in $K$ for every element of $\W$. The elements of
$\L$ which take $\mf C$ into itself will be called dominant and
the set of dominant elements will be denoted by $\L^+$. The
elements of the set $\L^-:=-\L^+$ will be called anti--dominant.


\subsection{}\label{convention}
To the affine root system $\Phi_{\rm af}$ we associate the unique
finite root system $\RR\subset V$ and affine root system $R$ as in
Section \ref{affine} for which the configuration of vanishing
hyperplanes of $\Phi_{\rm af}$ and $R$ coincides. Moreover, $\RR$
is chosen to be reduced if and only if $\Phi$ is reduced. In this
way, the Weyl groups associated to $\RR$ and $R$ are $\W=W(\Phi)$
and $W=W(\Phi_{\rm af})$ and we are able to freely use the
conventions and notations in Section \ref{affine}. The reflections
corresponding to the simple roots of $\Phi$, $\Phi_{\rm af}$ and
respectively $\RR$, $R$ will coincide although the roots
themselves may be different (but necessarily a scaling of each
other). Furthermore, the lattice $\L$ becomes an intermediary
lattice between $Q$ and $P$, the root lattice and, respectively,
the weight lattice associated to $\RR$.

\subsection{}
The length $\ell(w)$ of an element $w$ of $W^e$ is the number of
affine hyperplanes in $\A$  separating $\cal C$ and $w\cdot\cal
C$.  For example, the elements of $\Omega$ have length zero. For
each $w$ in $W$, $\ell(w)$ is the length of a reduced (i.e.
shortest) decomposition of $w$ in terms of the simple affine
reflections $r_0,\cdots, r_n$. The unique highest length element
of $\W$ is denoted by $w_{\circ}$ and, for $\l$ in $\L$ we denote
by $\w_\l$ the unique minimal length element of $\W$ for which
$\w_\l^{-1}\cdot \l$ is anti--dominant.We list below a few useful
 facts about the length function. For a proof see, for example,
 \cite[Lemma 1.1]{ion}.
\begin{Lm}\label{length}
Assume that $\l$ and $\mu$ are anti--dominant elements of $\L$ and
that $\w$ is an element of $\W$. Then the following are true,
\begin{enumerate}
\item $\ell(t_{\l+\mu})=\ell(t_\l)+\ell(t_\mu) $ \item
$\ell(t_\l\w)=\ell(t_\l)+\ell(\w)$
\end{enumerate}
\end{Lm}

 The elements of $\L$ which keep $x_0$ into the
chamber $\cal C$ are called minuscule and their set is denoted by
$\cal O_\L$. The set $\cal O_\L$ is finite, of the same
cardinality as the group $\Omega$. In fact, we can parameterize
$\Omega$ by the elements of $\cal O_\L$ as follows: for each
$\l\in \cal O_\L$ let $\omega_\l$ denote the unique element of
$\Omega$ for which $\omega_\l(0)=\l$. It is easy to see that
$\omega_\l=t_\l\w_\l$.

For each $\l$ in $\L$ let us denote by $\tilde\l$ the unique
element of $\cal O_\L$ which lies in the same orbit of the affine
action of $W$ on $\L$ (regarded as a subset of $\A$ via $\l\mapsto
\l+x_0$) and by $w_\l$ the unique minimal length element of $W$
for which $w_\l\cdot\tilde \l=\l$.

\subsection{} The groups $G$, $K$, $M$ and
$U$ are unimodular. We choose correspondingly the Haar measures
$dg$, $dk$, $dm$ and $du$ normalized such that they give volume
one to the intersection between the corresponding group and the
Iwahori $I$. The parabolic subgroup $P=MU$ is not unimodular; a
left invariant Haar measure on $P$ is given by
$$
\int_P f(p)d_\ell p :=\int_M\int_U f(mu)dmdu
$$
and a right invariant Haar measure is given by $d_r p:=\d_P(p)
d_\ell p$, where $\d_P:P\to \Re^\times_+$ is the modular function
of $P$. In fact, $\d_P(mu)=|{\rm det}(Ad_{|U}(m))|$ for any $m$ in
$M$ and $u$ in $U$; we denoted by $Ad_{|U}$ the automorphism of
the Lie algebra of $U$ given by the adjoint representation. Since
$\d_P$ is trivial on $U$ and $\ring {}M$, it can be seen as an
unramified character of $M$. For example, we will sometimes write
$\d_P(\mf e^\l)$ to refer to $\d_P(t_\l)$. With the above notation
the following formula holds
\begin{equation}\label{integration}
\int_G f(g)dg=\int_K\int_P f(pk)dkd_\ell p = \int_K\int_M\int_U
f(muk)dkdmdu
\end{equation}

For any element $w$ of $W^e$ the volume of the double coset $IwI$
is denoted by $\mf t(w)$. If $w$ has the reduced decomposition
$\omega r_{i_1}\cdots r_{i_l}$ with $\omega$ in $\Omega$ and
$r_{i_j}$ simple reflections, then $\mf t(w)=\mf t_{i_1}\cdots \mf
t_{i_l}$. The connection between these numbers and the values of
the modular function is the following (see \cite[Corollary
3.2.13]{macdonald}).
\begin{Lm}\label{lemma2}
Let $\l$ be an element of $\L^+$ and let $\w$ be in $\W$. Then,
\begin{equation}\label{eq1}
\d_P(t_\l)=\mf t(\w t_\l\w^{-1})
\end{equation}
\end{Lm}

\subsection{} We will list here the main
results about the structure of $G$ which will be used in this
paper. All the results are due to Bruhat and Tits
\cite{bruhat-tits} (see also \cite{cartier}, \cite{macdonald},
\cite{tits} for streamlined expositions).

\begin{Thm}[Iwasawa decomposition]\label{iwasawa}
We have $G=PK$. Moreover, we have the following disjoint unions
\begin{equation}\label{iwasawa1}
G=\bigcup_{\w\in \W} P \w I=\bigcup_{w\in W^e} \ring {}MUw I
\end{equation}
\end{Thm}
For second equality in (\ref{iwasawa1}) follows from the first
one, by keeping in mind that $P=MU$, $M/~\ring {}M\cong \L$ and
$W^e=\W\ltimes \L$.


\begin{Thm}[Iwahori decomposition] Let $I^0=I\cap M$,
$I^+=I\cap U$ and $I^-=I\cap U^-$. Then,
\begin{equation}\label{iwahori1}
I=I^0I^+I^-
\end{equation}
and the same holds with the product taken in any order. Moreover,
$I^0=\ring {}M$, $I^+=K\cap U$ and
\begin{equation}\label{iwahori2}
t_\l I^+ t^{-1}_\l \supset I^+ \quad \text{and} \quad t_\l I^-
t^{-1}_\l \subset I^-
\end{equation}
for any $\l\in \L^+$. Furthermore,
\begin{equation}\label{iwahori3}
U=\bigcup_{\l\in \L^+}t_\l I^+t_\l^{-1}
\end{equation}
\end{Thm}

\begin{Thm}[Bruhat decomposition] We have $G=INI$. More precisely, we
have the following disjoint union
\begin{equation}\label{bruhat1}
G=\bigcup_{w\in W^e} I w I
\end{equation}
The double cosets $IwI$ are called Bruhat cells and they have the
following compatibility with the group structure of $G$
\begin{eqnarray}
Ir_iI\cdot IwI\subset Ir_iwI\cup I w I & & \text{for any $w\in
W^e$
and $0\leq i\leq n$} \label{bruhat4} \\
Ir_iI\cdot IwI=Ir_iwI & & \text{if and only if
$\ell(r_iw)=\ell(w)+1$} \label{bruhat5}
\end{eqnarray}
\end{Thm}

Let us note that we also have $K=I(N\cap K)I$ and the following
disjoint union
\begin{equation}\label{bruhat2}
K=\bigcup_{\w\in \W} I^+ \w I
\end{equation}
Furthermore, $G$ can be written as the  disjoin union
\begin{equation}\label{bruhat3}
G=\bigcup_{\l\in \L} I t_\l K
\end{equation}

\subsection{} The extended affine Hecke
algebra associated to the pair $(G,I)$ is defined as $\mf
H^e:=C_c(I\backslash  G/I)$, the set of compactly supported
$I$--bi--invariant functions on $G$, endowed with the associative
algebra structure given by convolution
$$
(f*h)(x):=\int_G f(g)h(g^{-1}x)dg
$$
If we let $\T_w$ denote $1_{IwI}$, the characteristic function of
the double coset $IwI$ then, from the Bruhat decomposition,
$\{\T_w\}_{w\in W^e}$ is a basis of $\mf H^e$. To avoid cumbersome
notation we use the symbol $\T_i$ to refer to $\T_{r_i}$ for any
affine simple reflection $r_i$.

 We define the affine Hecke algebra $\mf H$, respectively the
finite Hecke algebra $\ring {\mf H}$, as being the subalgebra of
$\mf H^e$ with basis $\{\T_w\}_{w\in W}$, respectively with basis
$\{\T_{\w}\}_{\w\in \W}$. Again, from the Bruhat decomposition, we
see that $\ring {\mf H}=C_c(I\backslash  K/I)$. The characteristic
function of $I$ is the unit in all these algebras and it will be
denoted by $1$. The structure of the algebra $\mf H^e$ has been
described by Iwahori and Matsumoto.

\begin{Thm} With the above notation the following statements hold
\begin{enumerate}
\item The affine Hecke algebra $\mf H$ is generated by the
elements $\T_0, \T_1,\cdots, \T_n$ and the following relations are
a complete set of relations among them
\begin{enumerate}
\item the elements $\T_0, \T_1,\cdots, \T_n$ satisfy the same
braid relations as the simple affine reflections $r_0, r_1,\cdots,
r_n$;

\item for all $0\leq i\leq n$ we have
$$
\T_i^2=(\mf t_i-1)\T_i+ \mf t_i
$$
\end{enumerate}
\item For elements $w_1$, $w_2$ of $W^e$ such that
$\ell(w_1w_2)=\ell(w_1)+\ell(w_2)$ we have
$$
\T_{w_1w_2}=\T_{w_1}\T_{w_2}
$$
In particular, the elements $\{\T_\omega\}_{\omega\in \Omega}$
form a group isomorphic to $\Omega$ (and denoted by the same
symbol) and $\mf H^e=\Omega\ltimes \mf H$.
\end{enumerate}
\end{Thm}

The extended affine Hecke algebra admits a second description due
to Bernstein (unpublished) and Lusztig \cite{lusztig}. Before
stating the result we need the following notation. For $\l$ in
$\L$ we can find two elements $\mu_1$, $\mu_2$ in $\L^+$ such that
$\l=\mu_1-\mu_2$. Let $\X_\l$ denote the element $\mf
t(t_{\mu_1})^{-\frac{1}{2}}\mf
t(t_{\mu_2})^{\frac{1}{2}}\T_{t_{\mu_1}}\T^{-1}_{t_{\mu_2}}$ of
the extended affine Hecke algebra $\mf H^e$; $\X_\l$ does not
depend on the choice of $\mu_1$ and $\mu_2$. Moreover, the
elements $\{\X_\l\}_{\l\in \L}$ generate a commutative subalgebra
of $\mf H^e$ isomorphic to $\mf R$ and which will be denoted by
$\mf R_{\X}$.
\begin{Thm}\label{lusztigpresentation}
The extended affine Hecke algebra $\mf H^e$ is generated by the
finite Hecke algebra $\ring{\mf H}$ and the commutative algebra
$\mf R_{\X}$ such that the following relations are satisfied for
any $\l$ in $\L$ and any $1\leq i\leq n$
$$\hspace*{-3.7cm}
\T_i \X_\l - \X_{r_i(\l)}\T_i =(\mf
t_i-1)\frac{\X_\l-\X_{r_i(\l)}}{1-\X_{-\a_i}}~~ \text{ if }~
2\a_i\not \in R
$$
$$ \T_n \X_\l - \X_{r_n(\l)}\T_n  =\left(\mf
t_n-1+\mf t_n^{\frac{1}{2}}(\mf t_0^{\frac{1}{2}}-\mf
t_0^{-\frac{1}{2}})\X_{-\a_n}\right)\frac{\X_\l-\X_{r_n(\l)}}{1-\X_{-2\a_n}}~
~\text{ if }~~~ 2\a_n\in R
$$
\end{Thm}

\section{Nonsymmetric Macdonald polynomials}\label{nonsymm}

\subsection{} In this Section we will list the facts about nonsymmetric
Macdonald polynomials which are relevant to our problem. For a
detailed account of their construction and basic properties we
refer the reader to \cite{macbook} and \cite{c}. The results
concerning the limit as $q\to \infty$ of nonsymmetric Macdonald
polynomials were obtained in \cite{ion}.

 For a root system $\RR$ as in Section 1, the
nonsymmetric Macdonald polynomials $E_\l(q,t)$ are family of
orthogonal polynomials which form a basis of $\cal R_{q,t}$. These
polynomials are indexed by the weight lattice of $\RR$ and they
depend rationally on parameters $q$ and $t=(t_\a)_{\a\in R}$. They
were defined by Opdam \cite{opdam} and Macdonald for various
specializations of the parameters and in full generality by
Cherednik \cite{c}.  It was proved in \cite[Section 3.1]{ion} that
the coefficients of nonsymmetric Macdonald polynomials have finite
limits as $q$ approaches infinity. We will use the notation
$E_\l=E_\l(\infty,t)$ for the nonsymmetric Macdonald polynomials
in this limit.

To avoid introducing more notation we will state some properties
of the polynomials $E_\l$ using the notation used in the present
context. For this, we adhere to the convention in Section
\ref{convention} and, moreover, we restrict our attention to the
subring $\cal R_{t,\L}$ of $\cal R_t$ spanned by
$\{e^\l\}_{\l\in\L}$. The polynomials $\{E_\l\}_{\l\in \L}$ form a
basis of $\cal R_{t,\L}$ and the extended affine Hecke algebra
$\mf H^e$ acts on $\cal R_{t,\L}$ as follows
$$\hspace*{-1.6cm}
\T_i\cdot  e^\l = t_i  e^{r_i(\l)}+(t_i-1)\frac{ e^\l-
e^{r_i(\l)}}{1- e_{-\a_i}}~~ \text{ if }~ 1\leq i\leq n \text{ and
} 2\a_i\not \in R
$$
$$ \T_n\cdot  e^\l = t_n  e^{r_n(\l)}+\left(
t_n-1+ t_n^{\frac{1}{2}}( t_{03}^{\frac{1}{2}}-
t_{03}^{-\frac{1}{2}}) e^{-\a_n}\right)\frac{ e^\l-e^{r_n(\l)}}{1-
e^{-2\a_n}}~ ~\text{ if }~~~ 2\a_n\in R
$$
$$\hspace*{-7.1cm}
\X_\mu\cdot  e^\l=  e^{\l+\mu} ~~ \text{ for all }~ \mu \in \L
$$

If $\l$ is in $\L$ we denote by ${\bf t}^{(\l,\overline\l)}$ the
element of $\F$
$$
(t_{n}^*t_n)^{\frac{1}{2}(\l_-, \l_n^\vee)}
\prod_{i=1}^{n-1}t_i^{(\l_-, \l_i^\vee)}
$$
where $t_n^*$ equals $t_n$ if $R$ is reduced or $t_{01}$ if $R$ is
nonreduced. If $\mf A$ is the  affine Artin group of which the
affine Hecke algebra $\mf H$ is a quotient,  let us consider the
group morphism defined by
\begin{equation*}
\xi:\mf A\to \F_t, \quad \xi(\T_i)=t_i,\ i\not=0\ \ \text{and}\ \
\xi( \T_{0})=t_{01}^{\frac{1}{2}}t_{03}^{\frac{1}{2}}
\end{equation*}
We will abuse notation and write $\xi(w)$ to refer to $\xi(\T_w)$
for $w$ in $W$. Given  $\l$ in $\L$ let us define the following
normalization factor
\begin{equation}
j_\l:=\xi(w_\l){\bf t}^{(\l,\overline\l)}
\end{equation}
The following result was proved in \cite[Corollary 3.3]{ion}.
\begin{Prop}\label{recursion}
Let $\l$ be a in $\L$. Then,
$$
\T_{w_\l}\cdot {j_{\tilde \l}}e^{\tilde \l}={j_\l} E_\l
$$
\end{Prop}
We should also note that as a consequence of Corollary 3.7 in
\cite{ion} the polynomials $E_\l$ are free of the parameter
$t_{02}$.
\section{Unramified principal series}

\subsection{} Let $\chi$ be any unramified character of $M$. The
{\sl  unramified principal series representation} of $G$ indexed
by $\chi$, denoted by $(\nu_\chi,\mf U(\chi))$, is the
representation obtained by normalized induction from the
one--dimensional representation of $P$ given by $\chi$ (which is
extended to $P$ by letting $U$ act trivially). More precisely,
$$\mf U(\chi)=
\left\{ f:G\to \C ~ \vrule
\begin{array}{c}
(i) ~f(mug)=\d_P^{\frac{1}{2}}(m)\chi(m)f(g)~ ~ \text{for } m\in
M, ~u\in U, ~g\in G
\\
\hspace*{-5.5cm} (ii) ~ f \text{ locally constant}
\end{array}\right\}
$$
and $\nu_\chi$ is the right $G$--action by right translations. The
following result is due to Casselman (see also \cite{borel}).
\begin{Prop}
Let $\chi$ be an unramified character of $M$. Then
\begin{enumerate}
\item $(\nu_\chi,\mf U(\chi))$ is admissible and the subspace of
its Iwahori fixed vectors $\mf U(\chi)^I$ has dimension $|\W|$;

\item $\mf U(\chi^{-1})$  is isomorphic to the contragradient
representation $\widetilde{ \mf U(\chi)}$;

\item If $(\pi,V)$ is an irreducible, admissible representation of
$G$ such that the space of its Iwahori fixed vectors $V^I$ is
nontrivial, then there exists a $G$--embedding of $V$ into some
unramified principal series representation. Conversely, if $V$ is
a any nontrivial subrepresentation of $\mf U(\chi)$, then $V^I$ is
also nontrivial.
\end{enumerate}
\end{Prop}

\subsection{} It will be more convenient to study the space of {\sl universal
unramified principal series} $(\rho,\mf U)$ which is defined by
$$
\mf U:=\left\{ f:G\to \R ~ \vrule
\begin{array}{c}
(i) ~f(mug)=\d^{\frac{1}{2}}_P(m){\mf e}^{-v_M(m)}f(g)~ ~
\text{for } m\in M, ~u\in U, ~g\in G
\\
\hspace*{-5.7cm} (ii) ~ f \text{ locally constant}
\end{array}\right\}
$$
The group $G$ acts by right translations making $\mf U$ a smooth
right $G$--module. Also, $\mf U$ admits a $\R$--module structure
defined as follows
$$
(\mf e^\l\cdot f)(g):=\mf e^{w_\circ(\l)}f(g)
$$
where $f$ is an element of $\mf U$ , $\l$ in $\L$, $g$ in $G$. For
$\chi$ an unramified character of $M$, it is easy to see that $\mf
U\otimes_{\mf R}\C_{\chi^{-1}}$ is precisely  $\mf U(\chi)$, the
unramified principal series associated to $\chi$.

\subsection{}

The universal unramified principal series admits a second
construction as compact induction {\rm c-ind}$_{\hspace{2pt}\ring
{}MU}^{G}(\C)$ from the trivial representation. Our presentation
here follows \cite{chriss}, \cite{kott} where the arguments for
split groups were presented (although they work in complete
generality). Let us denote by $\mf V$ the space {\rm
c-ind}$_{\hspace{2pt}\ring {}MU}^{G}(\C)$, which is by definition
$$\mf V:=
\left\{ f:G\to \C ~ \vrule
\begin{array}{c}
\hspace{-1.1cm}(i) ~f(mug)=f(g)~ ~ \text{for } m\in \ring {}M,
~u\in U, ~g\in G
\\
\hspace{-4.8cm} (ii) ~ f \text{ locally constant}\\
(iii)~ \text{supp}(f) \subset \ring {}MUC_f ~ \text{ for some
compact set } C_f \text{ of } G
\end{array}\right\}
$$
It is easy to see that, in fact, $\mf
V=C_c^\infty(\hspace{2pt}\ring {}MU\backslash G)$. The action of
$G$ by right translations endows $\mf V$ with a smooth right
$G$--module structure. It can be also endowed with a left $\mf
R$--module structure (which commutes with the right $G$--action)
as follows. For $\l$ in $\L$, $f$ in $\mf V$ and any $g$ in $G$
define
\begin{equation}\label{RonV}
(\mf e^\l\cdot f)(h):=\d^{\frac{1}{2}}_P(\mf
e^{w_\circ(\l)})f(t_{w_\circ(\l)}^{-1}g)
\end{equation}

For a map $f:G\to  \mf R$ let us denote by $f_\l(\cdot): G\to \C$
the functions for which
$$
f(g)=\sum_{\l\in \L} f_\l(g) \mf e^{\l}
$$
for all $g$ in $G$. Let,
\begin{eqnarray*}
\zeta: \mf U\to \mf V~ , & \hspace{-2.8cm} \zeta(f)(g):=f_0(g)
\quad
\text{for all $f$ in } \mf U \text{ and } g~ \text{in}~ G\\
\eta: \mf V\to \mf U~ , & \eta(h)(g):=\sum_{\l\in
\L}\d^{-\frac{1}{2}}_P(\mf e^{\l})h(t_\l g)\mf e^{\l} \quad
\text{for all $h$ in } \mf V \text{ and } g~ \text{in}~ G
\end{eqnarray*}

\begin{Lm}
With the above notation, the maps $\zeta$ and $\eta$ are well
defined, right $G$--module, left $\mf R$--module homomorphisms.
Moreover, they are inverse to each other and therefore $\mf U$ and
$\mf V$ are isomorphic as $(\mf R,G)$--modules.
\end{Lm}
\begin{proof}
It is a straightforward check that $\mf U$ and $\mf V$ are indeed
$(\mf R,G)$--modules and the above maps are $(\mf R,G)$--linear
and inverse to each other.
\end{proof}

\subsection{} The subspace of Iwahori fixed vectors of the
universal unramified principal series is isomorphic to $\mf
V^I=C_c(\hspace{2pt}\ring {}MU\backslash G/I)$. Theorem
\ref{iwasawa} allows us to consider the natural basis for $\mf V
^I$ given by the characteristic functions $\{\mf
v_x:=1_{\hspace{2pt}\ring {}MU  xI}\}_{x\in W^e}$.

Since $\mf V^I$ is a left $\mf R$--submodule of $\mf V$. The $\mf
R$--action on the elements of the above basis is given by
\begin{equation}\label{RonVI}
\mf e^\l\cdot \mf v_x=\d^{\frac{1}{2}}_P(\mf e^{w_\circ(\l)})\mf
v_{t_{w_\circ(\l)} x}
\end{equation}
for any $\l$ in $\L$ and $x$ in $W^e$. The right $G$--action on
$\mf V$ gives rise to a right action, by convolution, of $\mf H^e$
on $\mf V^I$. For example, for $x$ and $w$ in $W^e$ we have
\begin{equation}\label{eq9}
(\mf v_{1}* \T_w)(x) = \int_G 1_{\hspace{2pt}\ring
{}MUI}(g)1_{IwI}(g^{-1}x) dg
\end{equation}
which is easily seen to be the volume of the set $UI\cap
xIw^{-1}I$. In particular, $c_{x,w}:=(\mf v_{1}* \T_w)(x)$ it is
always a non--negative integer.
\begin{Lm}\label{triangular}
With the above notation we have
$$
\mf v_{1}* \T_w=\sum_{x\leq w} c_{x,w}\mf v_x
$$
for some  non--negative integers $c_{x,w}$. Moreover, $c_{w,w}$ is
non-zero.
\end{Lm}
\begin{proof}
Since $\{\mf v_x\}_{x\in W^e}$ form a basis for $\mf V^I$ we can
derive the expansion of $\mf v_1* \T_w$ in this basis simply by
evaluating it at various elements of the extended Weyl group
$W^e$. Our claim follows if we show that $UI\cap xIw^{-1}I$ is
empty unless $x\leq w$.

Indeed, assume that $u$ is an element of $U$ such that $ux\in
IwI$. By (\ref{iwahori3}) we know that there exists $\l$ in $\L^-$
for which $t_\l ut_\l^{-1}$ lies in $I^+$. In such a case, $t_\l
ut_\l^{-1}\cdot t_\l x\in t_\l IwI$ and we can deduce that
\begin{equation}\label{eq2}
It_\l x I\subset It_\l IwI \subset \bigcup_{y\leq w} It_\l yI
\end{equation}
 The
last inclusion is justified by the multiplicative properties of
the Bruhat cells. In conclusion we must have $x\leq w $. The
coefficient $c_{x,w}$ is the volume of $UI\cap xIw^{-1}I$ which is
a non--negative integer. Furthermore, $c_{w,w}$ is strictly
positive since $I\subset UI\cap wIw^{-1}I$.
\end{proof}
In some cases we can obtain a sharper result.
\begin{Lm}\label{lemma1}
Let  $\mu$ be an element of $\L^-$ and $\w$ be an element of $\W$.
Then $$\mf v_1* \T_{t_\mu \w}=\mf v_{t_\mu\w}$$
\end{Lm}
\begin{proof}  If we proceed as in the above proof we
find that the relation (\ref{eq2}) becomes
$$
It_\l x I\subset It_\l It_\mu \w I = It_{\l+\mu} \w  I
$$
where the equalities  are  consequences of Lemma \ref{length} and
of the multiplicative properties of the Bruhat cells. Therefore,
$$
\mf v_1* \T_{t_\mu \w}=c_{t_\mu\w,t_\mu\w}\mf v_{t_\mu\w}
$$
Let us argue first that $$\mf v_1* \T_{t_\mu}=\mf v_{t_\mu}$$
Since, by the above argument,
$$
\mf v_1* \T_{t_\mu}=c_{t_\mu,t_\mu}\mf v_{t_\mu}
$$
we only need to show that $c_{t_\mu,t_\mu}=1$. Using the Iwahori
decomposition (\ref{iwahori1}), (\ref{iwahori2}) we obtain that
$$I\subset UI\cap t_\mu I t_\mu^{-1}I=UI\cap t_\mu I^- t_\mu^{-1}I\subset UI\cap U^-I =I$$
and we deduce that $c_{t_\mu,t_\mu}=1$.

Now, our hypothesis implies that $\T_{t_\mu
\w}=\T_{t_\mu}\T_{\w}$. Hence,
$$
\mf v_1* \T_{t_\mu \w}=\mf v_1* \T_{t_\mu}*\T_{ \w}=\mf v_{t_\mu}*
\T_{\w}
$$
We deduce that $c_{t_\mu\w,t_\mu\w}$ equals the volume of the set
$t_\mu^{-1}Ut_\mu I\cap \w I\w^{-1}I$. Keeping in mind that $M$
normalizes $U$ and that $\w$ and $I$ are included  in $K$ we
obtain
$$I\subset t_\mu^{-1} U t_\mu I\cap \w I\w^{-1} I\subset UI\cap K=I$$
which shows that $c_{t_\mu\w,t_\mu\w}=1$.
\end{proof}
By a similar argument we can prove the following.
\begin{Lm}
Let $1_K$ be the characteristic function of $K$.
\begin{enumerate}
\item For any  $\w$ in $\W$ we have
\begin{equation}\label{Tw1K}
\T_{\w} *1_K= \mf t(\w) 1_K
\end{equation}
\item For any $\nu$ in $\cal O_\L$ we have
\begin{equation}\label{TlK}
\T_{t_\nu} *1_K=  \mf t(\w_\nu)1_{It_\nu K}
\end{equation}
\end{enumerate}
\end{Lm}
\begin{proof}
For the first claim, note that $\T_{\w} *1_K$ is an element of
$C_c(I\backslash G/ K)$ and by  (\ref{bruhat3}) it is determined
by its values at the elements $t_\mu$. Using the convolution
formula we find that $(\T_{\w} *1_K)(t_\mu)$ equals the volume of
the set $I\w I\cap t_\mu K$. Recall that the representatives of
the elements of $\W$ are chosen to be in $K$ and therefore
$$
I\w I\cap t_\mu K\subset K\cap t_\mu K
$$
showing that $I\w I\cap t_\mu K$ is empty unless $\mu =0$, in
which case $I\w I\cap K=I\w I$. This set has volume $\mf t(\w)$
and our first statement is proved.

We argue similarly to prove the second statement. If $\mu$ is an
element of $\L$, the value of the integral
$$
(\T_{t_\nu} *1_K)(t_\mu)=\int_G 1_{It_\nu
I}(g)1_{K}(g^{-1}t_\mu)dg
$$
is zero, unless $t_\mu$ belongs to $It_\nu K$ or, equivalently,
unless $\mu=\nu$. When $\mu=\nu$ the value of the above integral
equals the volume of the set $t_\nu^{-1}It_\nu I\cap  K$. Since
$\nu$ is minuscule we know that $\omega_\nu=t_\nu\w_\nu$
normalizes $I$. Therefore, the volume of $t_\nu^{-1}It_\nu I\cap
K$ equals the volume of
$$\omega_\nu^{-1}I\omega_\nu
\w_\nu^{-1}I\cap K = I \w_\nu^{-1}I\cap K $$ which is precisely
$\mf t(\w_\nu)$.
\end{proof}
\subsection{} The following result reveals the close relationship
between $\mf H^e$ and $\mf V^I$. Again, the argument was presented
in \cite{chriss} and \cite{kott} (for split groups).
\begin{Prop}\label{HVI}
The right $\mf H^e$--linear map $$\varphi: \mf H^e\to \mf
V^I,\quad \varphi(\T_w):=\mf v_1* \T_w$$ is an isomorphism of
right $\mf H^e$--modules, the right $\mf H^e$--module $\mf V^I$ is
free of rank one and  the map $$\psi: \mf H^e\to {\rm End}_{\mf
H^e}(\mf V^I), \quad \psi(\T_w)(\mf v_1):=\mf v_1* \T_w
$$
is  an algebra isomorphism.
\end{Prop}
\begin{proof}
To justify the first claim, note that according to Lemma
\ref{triangular} the map $\varphi$ is lower triangular with
respect to the bases $\{\T_w\}_{w\in W^e}$ and $\{\mf v_w\}_{w\in
W^e}$ and therefore an an isomorphism of right $\mf H^e$--modules.
In consequence, $\mf V^I$ is a free right $\mf H^e$--module of
rank one and in consequence any $\mf H^e$--linear endomorphism of
$\mf V^I$ is uniquely determined by its value on $\mf v_1$.
Moreover, any element of ${\rm End}_{\mf H^e}(\mf V^I)$ is of the
form $\psi(h)$ for some $h\in \mf H^e$. Hence the map $\psi$ is
well defined a simple check confirms that it is an algebra
isomorphism.
\end{proof}
We will essentially restate this result in a form which will be
more suitable for our purposes.
\begin{Cor}
The right $\mf H^e$--linear map $$\varphi_\circ: \mf H^e\to \mf
V^I,\quad \varphi_\circ(\T_w):=\mf v_{w_\circ}* \T_w$$ is an
isomorphism of right $\mf H^e$--modules, and  the map
$$\psi_\circ: \mf H^e\to {\rm End}_{\mf H^e}(\mf V^I), \quad
\psi_\circ(\T_w)(\mf v_{w_\circ}):=\mf v_{w_\circ}* \T_w
$$
is  an algebra isomorphism.
\end{Cor}
\begin{proof} From Lemma \ref{lemma1} we know that $\mf v_{w_\circ}=\mf
v_1*\T_{w_\circ}$. The result now follows from the above
Proposition and the fact that $\T_{w_\circ}$ is invertible.
\end{proof}

The above Corollary has the following important consequence.
 Recall that $\mf V$ is also endowed  with a left $\mf R$--module
 which commutes with the right $\mf H^e$--action. Therefore, for any $\mu$ in $\L$, the
 left  multiplication by $\mf e^\mu$ gives a right $\mf
 H^e$--linear map
 $$
\gamma_\mu: \mf V^I\to \mf V^I, \quad \gamma_\mu(\mf v_x):=
\d_P^{\frac{1}{2}}(\mf e^{w_\circ(\mu)})\mf v_{t_{w_\circ(\mu)} x}
 $$
which, by the above Corollary, must be the image by $\psi_\circ$
of an  unique element of $\mf H^e$, denoted by $\gamma(\mf
e^\mu)$. In this way we define an algebra map
$$
\gamma: \mf R\to \mf H^e
$$
\begin{Prop}\label{XR} Let  $\mu$ be any element of $\L$.
Then, with the above notation,
$$\gamma(\mf e^\mu)=\X_\mu$$
\end{Prop}
\begin{proof} Let us remark first that, since $\gamma$ is an
algebra isomorphism we only need to check our claim for elements
of $\L^+$. Assume now that $\mu$ is dominant. We will show that
$$\mf v_{w_\circ}*\X_\mu=\gamma_\mu(\mf v_{w_\circ})=\d_P^{\frac{1}{2}}(\mf
e^{w_\circ(\mu)})\mf v_{w_\circ t_\mu }$$ which is enough to
conclude our proof.

Because $\mu$ is dominant we  have $\X_\mu=\mf
t(t_\mu)^{-\frac{1}{2}}\T_{t_\mu}$. Therefore, after taking  into
account Lemma \ref{lemma2}, our claim is reduced to showing that
$$
\mf v_{w_\circ}*\T_{t_\mu}=\mf v_{w_\circ t_\mu }
$$
But, $\mf v_{w_\circ}=\mf v_1*\T_{w_\circ}$ and
$\T_{w_\circ}\T_{t_\mu}=\T_{w_\circ t_\mu}=\T_{t_{w_\circ(\mu)}
w_\circ }$. The element $w_\circ(\mu)$ is in $\L^-$ and therefore
Lemma \ref{lemma1} applies allowing us reach the conclusion.
\end{proof}

\section{The Satake transform}

\subsection{} Let us consider  $\mf V^K=C_c(\hspace{2pt}\ring
{}MU\backslash G/K)$, the space of $K$--fixed vectors in $\mf V$.
A basis for $\mf V^K$ is given by the characteristic functions
$\{\mf u_\l:=1_{\hspace{2pt}\ring {}MUt_{w_\circ(\l)} K}\}_{\l\in
\L}$. This can be easily seen using the Iwasawa decomposition
$G=PK=MUK$ and the fact that $M/~\ring {}M\cong \L$. The left $\mf
R$--action on $\mf V$ transfers to a left $\mf R$--action on $\mf
V^K$
$$
\mf e^\mu \cdot \mf u_\l=\d^{\frac{1}{2}}_P(\mf
e^{w_\circ(\l)})\mf u_{\mu+\l}
$$
In particular,  the map
\begin{equation}\label{VKtoR}
\sigma:\mf V^K\to \mf R~, \quad \sigma(\mf
u_\l):=\d^{-\frac{1}{2}}_P(\mf e^{w_\circ(\l)})\mf e^\l
\end{equation}
 is an isomorphism of
left $\mf R$--modules. In other words, $\sigma(\mf u_\l)$ is the
unique element of $\mf R$ for which $\sigma(\mf u_\l)\cdot \mf
u_0=\mf u_\l$. The connection between the canonical bases of $\mf
V^I$ and $\mf V^K$ is the following.
\begin{Lm}\label{lemma3}
For any $\l$ in $\L$ and any $\w$ in $\W$ we have
$$
\mf  t(\w)\mf u_{w_\circ(\l)}=\mf v_{t_\l\w}* 1_K
$$
\end{Lm}
\begin{proof}
Obviously $\mf v_{t_\l\w}* 1_K$ is an element of $\mf V^K$ and the
coefficients of its expansion in terms of the basis of  $\mf V^K$
given above are obtained by evaluating $\mf v_{t_\l}* 1_K$ at the
elements $t_\mu$. We easily find that $(\mf v_{t_\l}* 1_K)(t_\mu)$
is the volume of the set $\ring {}MUt_\l \w I\cap t_\mu K$. Now,
$$
\ring {}MUt_\l \w I\cap t_\mu K\subset \ring {}MUt_\l K\cap \ring
{}MUt_\mu K
$$
and therefore $\ring {}MUt_\l I\cap t_\mu K$ is empty unless
$\mu=\l$, in which case $$\ring {}MUt_\l \w I\cap t_\l K=t_\l(U
\cap K)\w I=t_\l I\w I$$ For the last equality we used the fact
that $U$ is normalized by $M$ and also the equality $I\w I=I^+\w
I$. Since the volume of $I\w I$ equals $\mf t(\w)$ our claim is
proved.
\end{proof}
\begin{Lm} Let $\l$ be in $\L$. Then,
$$
\mf t(w_\circ)^{-1}\mf v_{w_\circ}*\X_\l*1_K=
\d^{\frac{1}{2}}_P(\mf e^{w_\circ(\l)})\mf u_\l
$$
\end{Lm}
\begin{proof} From Proposition \ref{XR} we know that
$\mf v_{w_\circ}*\X_\l=\d^{\frac{1}{2}}_P(\mf
e^{w_\circ{(\l)}})\mf v_{w_\circ t_{\l}}$. Now, Lemma \ref{lemma3}
gives us the desired result.
\end{proof}
\subsection{} The restriction of the map $\varphi_\circ$ from Proposition
\ref{HVI} to $C_c(I\backslash G/ K)$ takes values in $\mf V^K$. We
will denote this restriction by the same symbol.
\begin{Prop}
The map
$$
\Upsilon :C_c(I\backslash G/ K)\to \mf V^K, \quad f\mapsto \mf
t(w_\circ)^{-1}\mf v_{w_\circ}*f
$$
is a linear isomorphism.
\end{Prop}
\begin{proof}
Of course, the above map is injective since it is, modulo a
nonzero scalar factor, the restriction of an injective map. As for
the surjectivity, the above Lemma shows that for any $\l$ in $\L$
\begin{equation}\label{eq5}\Upsilon( \X_\l * 1_K)=\d^{\frac{1}{2}}_P(\mf e^{w_\circ(\l)})\mf
u_\l\end{equation} Therefore, all the elements of the basis of
$\mf V^K$ are in the image of $C_c(I\backslash G/ K)$.
\end{proof}

\begin{Def} The linear isomorphism
$$\Xi:=\sigma\circ \Upsilon :C_c(I\backslash G/K)\to \mf R$$
will be called the Satake transform.
\end{Def}

The above Proposition appeared previously in the literature \cite[Proposition 2.5]{kato}\footnote{I am grateful to Professor Shin-ichi Kato for pointing this out.}.
In \cite{kato},  $\Xi$ is called the Satake map.
The Satake transform can be equivalently described as follows: if
$f$ is an element  of $C_c(I\backslash G/K)$ then $\Xi(f)$  is the
unique element of $\mf R$ for which
\begin{equation}\label{eq7}
\mf t(w_\circ)^{-1}\mf v_{w_\circ}*f=\Xi(f) \cdot \mf u_0
\end{equation}

The terminology can be justified as follows. Let $f$ be an element
of $C_c(I\backslash G/K)$. Its Satake transform being an element
of $\mf R$, it can be written as
\begin{equation}\label{eq3}
\Xi(f)=\sum_{\l\in \L}c_{w_\circ(\l)} \mf e^{\l}
\end{equation}
With this notation $\Xi(f) \cdot \mf u_0=\sum_\l c_{w_\circ(\l)}
\d^{\frac{1}{2}}_P(\mf e^{w_\circ(\l)})\mf u_{\l}$ and
$$
(\Xi(f) \cdot \mf u_0)(t_{w_\circ(\mu)})=c_{w_\circ(\mu)}
\d^{\frac{1}{2}}_P(\mf e^{w_\circ(\mu)})
$$
On the other hand, using the formula (\ref{integration}) we obtain
$$
(\mf v_{w_\circ}*f)(t_{w_\circ(\mu)}) = \int_{U} f(w_\circ^{-1}
u^{-1}  t_{w_\circ(\mu)}) du
$$
and we can conclude that
\begin{equation}\label{eq4}
c_{w_\circ(\mu)}=\mf t(w_\circ)^{-1}\d^{-\frac{1}{2}}_P(\mf
e^{w_\circ(\mu)})\int_{U} f(w_\circ^{-1}ut_{w_\circ(\mu)}) du
\end{equation}
Examining (\ref{eq3}) and (\ref{eq4}) we can see that the
restriction of the map $\Xi$ to the spherical subalgebra
$C_c(K\backslash G/K)$ is (modulo the action of $w_\circ$)  the
classical Satake transform (see \cite{satake} or formula (3.3.4)
and Theorem 3.3.6 in \cite{macdonald}).

Let us indicate two immediate consequences of the above
considerations. First, if $f$ is a function in $C_c(I\backslash
G/K)$ which takes only non--negative values then all the
coefficients $c_\mu$ are non--negative. Second, the formula
(\ref{eq5}) shows that for all $\l$ in $\L$
\begin{equation}\label{eq6}
\Xi(\X_\l*1_K)=\mf e^\l
\end{equation}

\subsection{} The space $C_c(I\backslash G/K)$ is a left module for
extended affine Hecke algebra $\mf H^e$ which acts by convolution.
The Satake transform being a linear isomorphism, there is a unique
 left action of $\mf H^e$ on $\mf R$, denoted by $\cdot$~, which makes
 $\Xi$ left $\mf H^e$--linear.
A crude description of this action on $\mf R$ can be obtained
merely by  reexamining our earlier computations.
\begin{Lm}
Let $H$ be an element of $\mf H^e$ and let $\l$ be in $\L$. Then,
$$
H\cdot \mf e^\l=\Xi(H\X_\l *1_K)
$$
\end{Lm}
\begin{proof}
Straightforward from (\ref{eq6})  and the definition of the left
action of $\mf H^e$.
\end{proof}
However, we can obtain a detailed description of the action of the
generators of $\mf H^e$ in the Bernstein--Lusztig presentation.
\begin{Prop}
 The action of $\mf H^e$ on $\mf R$ is
completely described by the following formulas
$$\hspace*{-1.6cm}
\T_i\cdot \mf e^\l =\mf t_i \mf e^{r_i(\l)}+(\mf t_i-1)\frac{\mf
e^\l-\mf e^{r_i(\l)}}{1-\mf e_{-\a_i}}~~ \text{ if }~ 1\leq i\leq
n \text{ and } 2\a_i\not \in R
$$
$$ \T_n\cdot \mf e^\l =\mf t_n \mf e^{r_n(\l)}+\left(\mf
t_n-1+\mf t_n^{\frac{1}{2}}(\mf t_0^{\frac{1}{2}}-\mf
t_0^{-\frac{1}{2}})\mf e^{-\a_n}\right)\frac{\mf e^\l-\mf
e^{r_n(\l)}}{1-\mf e^{-2\a_n}}~ ~\text{ if }~~~ 2\a_n\in R
$$
$$\hspace*{-7cm}
\X_\mu\cdot \mf e^\l= \mf e^{\l+\mu} ~~ \text{ for all }~ \mu \in
\L
$$
\end{Prop}
\begin{proof}
By Theorem \ref{lusztigpresentation} the action of $\mf H^e$ on
$\mf R$ is completely determined if we describe how the elements
$\T_1, \cdots, \T_n$ and $\{\X_\mu\}_{\mu\in \L}$ act. We only
have to check that these elements have the predicted action. Let
us consider first the action of $\X_\mu$ for some $\mu$ in $\L$.
The above Lemma  and the formula (\ref{eq6}) gives us
\begin{eqnarray*}
\X_\mu\cdot \mf e^\l &=& \Xi({\X_{\mu+\l} *1_K}) \\
&=& \mf e^{\mu+\l}
\end{eqnarray*}
If we consider the action of $\T_i$ (for $i$ such that $2\a_i\not
\in R$, for example) the same formulas lead to
\begin{eqnarray*}
\T_i\cdot \mf e^\l &=&  \Xi(\T_i\X_\l *1_K)\\
&=& \Xi\left(\X_{r_i(\l)}\T_i*1_K +(\mf
t_i-1)\frac{\X_\l-\X_{r_i(\l)}}{1-\X_{-\a_i}}*1_K\right)\\
&=& \mf t_i \mf e^{r_i(\l)}+(\mf t_i-1)\frac{\mf e^\l-\mf
e^{r_i(\l)}}{1-\mf e_{-\a_i}}
\end{eqnarray*}
which is precisely the predicted action.
\end{proof}

\subsection{} In the light of (\ref{bruhat3}), a natural basis for the space
$C_c(I\backslash G/K)$ is given by the characteristic functions
$\{1_{It_\l K}\}_{\l\in \L}$.  For our purposes it will be more
convenient to choose the representative $w_\l\omega_{\tilde\l}$
for the coset $It_\l K$. This is possible since $It_\l
K=Iw_\l\omega_{\tilde\l} K$ for any $\l$ in $\L$. We will study
next the Satake transforms $\cal E_\l:=\Xi(1_{Iw_\l
\omega_{\tilde\l}K})$ which give rise to a natural basis of $\mf
R$.

The coefficients of $\cal E_\l$ are given by the formula
(\ref{eq4}). However, let us remark that in this situation one can
compute $(\mf v_{w_\circ}*1_{It_\l K})(t_{w_\circ(\mu)})$ in the
same way as for (\ref{eq9}) and obtain that it equals the volume
of the set $Uw_\circ I\cap t_{w_\circ(\mu)}Kt_{-\l} I$. The latter
has the same volume as $ Uw_\circ t_{-\mu}I\cap Kt_{-\l} I$ (use
the fact that $M$ normalizes $U$) or $ U^{-} t_{-\mu}I\cap
Kt_{-\l} I$ (use $w_\circ^{-1}Uw_\circ=U^-$ and $w_\circ\in K$)
and hence we can rewrite (\ref{eq4}) as
\begin{equation}\label{eq10}
c_{w_\circ(\mu)}=\mf t(w_\circ)^{-1}\d^{-\frac{1}{2}}_P(\mf
e^{w_\circ(\mu)})~ vol(U^- t_{-\mu}I\cap Kt_{-\l}I)
\end{equation}

\begin{Prop}\label{prop2}
Let $\nu$ be a minuscule element of $\L$. Then
$$
\cal E_\nu=\d_P^{-\frac{1}{2}}(t_\nu) \mf e^\nu
$$
\end{Prop}
\begin{proof}
 Since $\nu$ is minuscule, the
formula (\ref{TlK}) gives us
$$1_{It_\nu K}=\mf t(\w_\nu)^{-1} \T_{t_\nu}*1_K=\mf t(t_\nu)^{-\frac{1}{2}}\X_\nu*1_K$$
and (\ref{eq6}) implies the desired result. We have implicitly
used the equality $\mf t(t_\nu)=\mf t(\w_\nu)$ which follows from
the fact that $\omega_\nu=t_\nu \w_\nu$ normalizes $I$.
\end{proof}
Regarding the elements $w_\l$ we know that (see, for example,
\cite[Lemma 1.6]{ion}) that if $r_i\cdot \l\neq \l$ then
$w_{r_i\cdot \l}=r_iw_\l$.
\begin{Prop}\label{prop1}
Let $\l$ be an element of $\L$ and let $r_i$ be a simple affine
reflection for which $\ell(w_{r_i(\l)})>\ell(w_\l)$. Then,
$$
\T_i\cdot \cal E_\l= \cal E_{r_i(\l)}
$$
\end{Prop}
\begin{proof} Let us compute first $\T_i*1_{Iw_\l \omega_{\tilde\l}K}$. This is an
element of $C_c(I\backslash G/K)$ and hence it is completely
described by its values at the elements $w_\mu\omega_{\tilde\mu}$.
Now,
$$
(\T_i*1_{Iw_\l K})(w_\mu\omega_{\tilde\mu})=\int_G
1_{Ir_iI}(g)1_{Iw_\l K}(g^{-1}w_\mu\omega_{\tilde\mu})dg
$$
and the value of the integral is zero unless
$w_\mu\omega_{\tilde\mu}$ belongs to $Ir_iI\cdot Iw_\l
\omega_{\tilde\l}K$. Our hypothesis allows us to use
(\ref{bruhat5}) which implies that the product of the sets $Ir_i
I$ and $Iw_\l I$ equals $Iw_{r_i(\l)}I$. Hence, the value of the
above integral equals zero unless $\mu=r_i(\l)$. The value of the
integral at $w_{r_i(\l)}$ is easily seen to be the volume of the
set
$$\omega_{\tilde\l}^{-1}w_{\l}^{-1}r_iIr_i I\cap K\omega_{\tilde\l}^{-1}
w^{-1}_\l I$$
Note that by (\ref{bruhat4}) we have $r_iIr_iI\subset Ir_iI\cdot
Ir_iI\subset I\cup Ir_iI$ which implies that
$$
w_{\l}^{-1}r_iIr_i I\subset w^{-1}_\l I\cup w^{-1}_\l Ir_iI
$$
We observe that $\omega_{\tilde\l}^{-1}w^{-1}_\l I\cap
K\omega_{\tilde\l}^{-1}w^{-1}_\l I=\omega_{\tilde\l}^{-1}w^{-1}_\l
I$ and
\begin{eqnarray*}
\omega_{\tilde\l}^{-1}w^{-1}_\l Ir_iI \cap K\omega_{\tilde\l}^{-1}w^{-1}_\l I &
\subset & K\omega_{\tilde\l}^{-1}w^{-1}_\l Ir_iI \cap K\omega_{\tilde\l}^{-1}w^{-1}_\l I\\
& = & K\omega_{\tilde\l}^{-1}w^{-1}_{r_i(\l)} I \cap
K\omega_{\tilde\l}^{-1}w^{-1}_\l I
\end{eqnarray*}
which is of course empty. We conclude that
$$
\omega_{\tilde\l}^{-1}w^{-1}_\l I \subset
\omega_{\tilde\l}^{-1}w_{\l}^{-1}r_iIr_i I\cap
K\omega_{\tilde\l}^{-1}w^{-1}_\l I \subset
\omega_{\tilde\l}^{-1}w^{-1}_\l I
$$
and the volume of the set
$\omega_{\tilde\l}^{-1}w_{\l}^{-1}r_iIr_i I\cap
K\omega_{\tilde\l}^{-1}w^{-1}_\l I$ is one. This shows that
$$
\T_i*1_{Iw_\l \omega_{\tilde\l}K}=1_{Iw_{r_i(\l)}\omega_{\tilde\l}
K}
$$

\noi The Satake transform is a $\mf H^e$-linear map, implying
\begin{eqnarray*}
\T_i\cdot \cal E_{\l} &=& \Xi(\T_i*1_{Iw_\l \omega_{\tilde\l}K})\\
&=& \Xi( 1_{Iw_{r_i(\l)}\omega_{\tilde\l} K})
\end{eqnarray*}
which is our desired result.
\end{proof}
\begin{Cor}\label{cor1}
Let $\l$ be an element of $\L$. Then,
$$
\T_{w_\l}\cdot  \d_P^{-\frac{1}{2}}(\mf e^{\tilde \l})\mf
e^{\tilde \l}=\cal E_\l
$$
\end{Cor}
\begin{proof}
Straightforward from Proposition \ref{prop1} and Proposition
\ref{prop2}.
\end{proof}
The next result describes the relationship between the elements
$\cal E_\l$ and the nonsymmetric Macdonald polynomials $E_\l$. Let
us consider first the $\C$--algebra morphism $\cal R_{t,\L}\to \R$
which sends each $e^\l$ to $\mf e^\l$, the parameters
$t_{01},~t_{02},~t_{03},~t_1,\cdots,t_n$ to the numbers $\mf
t_n^*,~\mf t_0^*,~\mf t_0,~ \mf t_1,\cdots,\mf t_n$ respectively.
 If $f$ is an element of
$\cal R_{t,\L}$ its image through the above morphism will be
denoted by $f(\mf t)$.
\begin{Thm}\label{thm2} Let $\l$ be an element of $\L$. Then
$$\cal E_\l=j_\l E_\l(\mf t)$$
\end{Thm}
\begin{proof}
Straightforward from Corollary \ref{cor1} and Proposition
\ref{recursion}.
\end{proof}

The above result has as a consequence a new formula for weight
multiplicities of Demazure modules associated to irreducible
representations of complex, connected, simple algebraic groups
\cite{ion3}.

\section{Matrix coefficients}

\subsection{}\label{coeff} Let us define the following vectors of the unramified principal series
representation $\mf U(\chi)$. First, let $v_K$ be the $K$--fixed
vector of $\mf U(\chi)$ which takes the value
$\d_P^{\frac{1}{2}}(\mf e^\l)\chi(\mf e^\l)$ on the coset $\ring
{}MUt_\l K$ for each $\l$ in $\L$. Second, let $v_I$ be the
$I$--fixed vector of $\mf U(\chi)$ which takes the value
$\d_P^{\frac{1}{2}}(\mf e^\l)\chi(\mf e^\l)$ on the coset $\ring
{}MUt_\l w_\circ I$ for each $\l$ in $\L$ and, it vanishes on the
cosets $\ring {}MUt_\l y I$ for each $\l$ in $\L$ and $y\neq
w_\circ$ in $\W$. It is a straightforward check that $v_K$ and
$v_I$ are the vectors in $\mf U(\chi)$ which correspond
respectively (via the map $\eta$ and after tensoring with
$\C_\chi^{-1}$ ) to the vectors $\mf u_0$ and $\mf v_{w_\circ}$ in
$\mf V$. With this notation, the formula (\ref{eq7}) reads as
follows. For $f$ in $C_c(I\backslash G/K)$ we have
\begin{equation}\label{eq8}
v_I*f=\chi^{-1}(\Xi(f))v_K
\end{equation}

\subsection{} The contragradient representation of $\mf U(\chi)$
is $\mf U(\chi^{-1})$. Therefore, there exists a canonical,
non--degenerate, $G$--equivariant pairing
$$
\<\cdot,\cdot\>~: \mf U(\chi^{-1})\otimes \mf U(\chi)\to \C
$$
The spaces of $K$--fixed vectors on $\mf U(\chi)$ and $\mf
U(\chi^{-1})$ are 1--dimensional. We denote by $\tilde v_K$ the
unique $K$--fixed vector in $\mf U(\chi^{-1})$ such that
$$
\<\tilde v_K,v_K\>=1
$$

Let us define the following matrix coefficient
$$
E_\chi:G\to \C,\quad E_\chi(g):=\<\tilde v_Kg,v_I\>
$$
Since $\tilde v_K$ is $K$--fixed and $v_I$ is $I$--fixed the
function $E_\chi$ belongs to $C_c(K\backslash G/I)$ and it is
therefore completely determined by its values at the elements
$t_\l$ for $\l$ in $\L$. We are now ready to state and prove our
main result.
\begin{Thm}\label{thm1}
Let $\chi$ be an unramified character and $\l$ and element of
$\L$. Then,

$$
E_\chi(t_{-\l})=\frac{j_{\l}(\mf t)}{vol(Kt_{-\l} I)}\chi^{-1}(
E_{\l}(\mf t))
$$\end{Thm}

\begin{proof}
By Theorem \ref{thm2}, it is enough to show that
$$
E_\chi(t_{-\l})=vol(Kt_{-\l} I)^{-1}\chi^{-1}(\cal E_{\l})
$$
Let us remark first that
\begin{eqnarray*}
(E_\chi*f)(g) &=& \int_G E_\chi(gx^{-1})f(x)dx \\
&=& \<\tilde v_Kg,v_I*f\>
\end{eqnarray*}
In particular, for $f\in C_c(I\backslash G/K)$ and $g=1$ we obtain
that
$$
(E_\chi*f)(1)=\chi^{-1}(\Xi(f))
$$
Furthermore, if $f=1_{It_{\l}K}$ the above formula becomes
$$
vol(Kt_{-\l} I)E_\chi(t_{-\l})=\chi^{-1}(\cal E_{\l})
$$
\end{proof}


\end{document}